\documentclass[12pt]{article}
\usepackage{amsmath,amssymb,geometry,verbatim,enumerate,float,tikz}
\usetikzlibrary{calc,positioning,arrows}
\newtheorem{thm}{Theorem}

\newtheorem{prop}[thm]{Proposition}
\newtheorem{claim}{Claim}

\usepackage{setspace}
\usepackage{lineno}

\begin{document}

%\linenumbers
\onehalfspace

\title{Forests and Trees among Gallai Graphs}

\author{Felix Joos$^1$, Van Bang Le$^2$, Dieter Rautenbach$^1$}

\date{}

\maketitle

\begin{center}
$^1$
Institut f\"{u}r Optimierung und Operations Research, Universit\"{a}t Ulm, Ulm, Germany\\
\texttt{felix.joos@uni-ulm.de}, \texttt{dieter.rautenbach@uni-ulm.de}\\[3mm]

Institut f\"{u}r Informatik, Universit\"{a}t Rostock, Rostock, Germany\\
\texttt{le@informatik.uni-rostock.de}
\end{center}

\begin{abstract}
The Gallai graph $\Gamma(G)$ of a graph $G$ has the edges of $G$ as its vertices
and two distinct vertices $e$ and $f$ of $\Gamma(G)$ are adjacent in $\Gamma(G)$
if the edges $e$ and $f$ of $G$ are adjacent in $G$ but do not span a triangle in $G$.
Clearly, $\Gamma(G)$ is a subgraph of the line graph of $G$.
While line graphs can be recognized efficiently 
the complexity of recognizing Gallai graphs is unknown.
In the present paper we characterize those graphs whose Gallai graphs are forests or trees, respectively.

\bigskip

\end{abstract}

\pagebreak

\section{Introduction}

We consider finite, simple, and undirected graphs and use standard terminology and notation \cite{west}.
For a graph $G$, the {\it Gallai graph} $\Gamma(G)$ of $G$ has the edges of $G$ as its vertices,
that is, $V(\Gamma(G))=E(G)$,
and two distinct vertices $e$ and $f$ of $\Gamma(G)$ are adjacent in $\Gamma(G)$
if the edges $e$ and $f$ of $G$ are adjacent in $G$ but do not span a triangle in $G$.
Gallai graphs were introduced by Gallai \cite{ga} in connection with cocomparability graphs
and were used by Chv\'{a}tal and Sbihi \cite{cs} in their polynomial time recognition algorithm 
for claw-free perfect graphs.
Obviously, 
the Gallai graph $\Gamma(G)$ is a spanning subgraph of the well-known {\it line graph} $L(G)$ of $G$ \cite{west}.
The {\it anti-Gallai graph} or {\it triangular line graph} $\Delta(G)$ of $G$
is the complement of $\Gamma(G)$ in $L(G)$,
that is, $V(\Delta(G))=E(G)$ and $E(\Delta(G))=E(L(G))\setminus E(\Gamma(G))$.
Anti-Gallai graph were introduced by Jarret \cite{jarrett}.

Gallai and anti-Gallai graphs were studied in \cite{lakshmanan,le1988,le1994,le}.
While the recognition of line graphs can be done efficiently \cite{lehot,roussopoulos},
it is hard to recognize anti-Gallai graphs \cite{anand-etal}
and the complexity of recognizing Gallai graphs is an open problem.
The characterizations of Gallai graphs given by Le \cite{le} 
do not seem to lead to an efficient recognition algorithm.
Therefore, 
further insight into the structure of Gallai graphs
and efficiently checkable characterizations of subclasses of Gallai graphs are of interest.
In the present paper we prove the following two results
characterizing those graphs whose Gallai graphs are forests or trees, respectively.

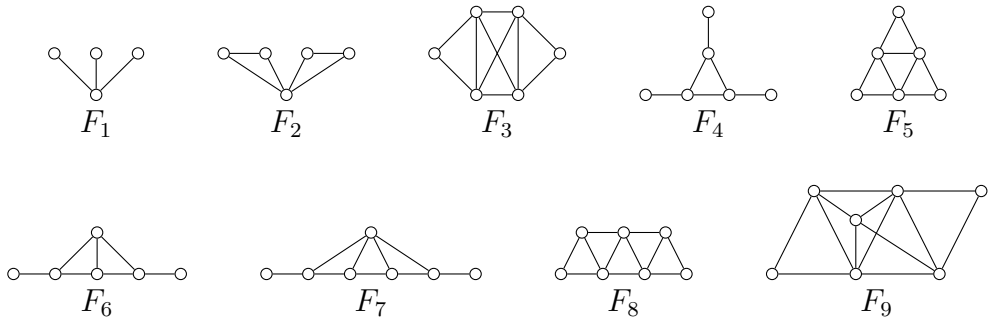
\begin{figure}[H]
\begin{center}
\begin{tikzpicture}[scale=.55] %,thick]
%\tikzstyle{vertex}=[circle,inner sep=2pt,fill=black];
\tikzstyle{vertex}=[draw,circle,inner sep=1.5pt] 
\node[vertex] (1) at (0,2)  {};
\node[vertex] (2) at (1,2)  {};
\node[vertex] (3) at (2,2)  {}; 
\node[vertex] (4) at (1,1)  {}; 

\draw (1) -- (4) -- (2);
\draw (3) -- (4);

\draw(1,.3) node {$F_1$};
\end{tikzpicture}
\qquad
\begin{tikzpicture}[scale=.55] %,thick]
%\tikzstyle{vertex}=[circle,inner sep=2pt,fill=black];
\tikzstyle{vertex}=[draw,circle,inner sep=1.5pt] 
\node[vertex] (1) at (0,2)  {};
\node[vertex] (2) at (1,2)  {};
\node[vertex] (3) at (2,2)  {}; 
\node[vertex] (4) at (3,2)  {}; 
\node[vertex] (5) at (1.5,1)  {}; 

\draw (1) -- (2) -- (5) -- (1);
\draw (3) -- (4) -- (5) -- (3);

\draw(1.5,.3) node {$F_2$};
\end{tikzpicture}
\qquad
\begin{tikzpicture}[scale=.55] %,thick]
%\tikzstyle{vertex}=[circle,inner sep=2pt,fill=black];
\tikzstyle{vertex}=[draw,circle,inner sep=1.5pt] 
\node[vertex] (1) at (0,2)  {};
\node[vertex] (2) at (1,1)  {};
\node[vertex] (3) at (2,1)  {}; 
\node[vertex] (4) at (3,2)  {}; 
\node[vertex] (5) at (2,3)  {}; 
\node[vertex] (6) at (1,3)  {}; 

\draw (1) -- (2) -- (3) -- (4) -- (5) -- (6) -- (1);
\draw (2) -- (5) -- (3) -- (6) -- (2);

\draw(1.5,.3) node {$F_3$};
\end{tikzpicture}
\qquad
\begin{tikzpicture}[scale=.55] %,thick]
%\tikzstyle{vertex}=[circle,inner sep=2pt,fill=black];
\tikzstyle{vertex}=[draw,circle,inner sep=1.5pt] 
\node[vertex] (1) at (0,1)  {};
\node[vertex] (2) at (1,1)  {};
\node[vertex] (3) at (2,1)  {}; 
\node[vertex] (4) at (3,1)  {}; 
\node[vertex] (5) at (1.5,2)  {}; 
\node[vertex] (6) at (1.5,3)  {}; 

\draw (1) -- (2) -- (3) -- (4);
\draw (2) -- (5) -- (3);
\draw (5) -- (6);

\draw(1.5,.3) node {$F_4$};
\end{tikzpicture}
\qquad
\begin{tikzpicture}[scale=.55] %,thick]
%\tikzstyle{vertex}=[circle,inner sep=2pt,fill=black];
\tikzstyle{vertex}=[draw,circle,inner sep=1.5pt] 
\node[vertex] (1) at (0,1)  {};
\node[vertex] (2) at (1,1)  {};
\node[vertex] (3) at (2,1)  {}; 
\node[vertex] (4) at (0.5,2)  {}; 
\node[vertex] (5) at (1.5,2)  {}; 
\node[vertex] (6) at (1,3)  {}; 

\draw (1) -- (2) -- (3) -- (5) -- (6) -- (4)-- (1);
\draw (2) -- (4) -- (5) -- (2);

\draw(1,.3) node {$F_5$};
\end{tikzpicture}
\end{center}

\begin{center}
\begin{tikzpicture}[scale=.55] %,thick]
%\tikzstyle{vertex}=[circle,inner sep=2pt,fill=black];
\tikzstyle{vertex}=[draw,circle,inner sep=1.5pt] 
\node[vertex] (1) at (0,1)  {};
\node[vertex] (2) at (1,1)  {};
\node[vertex] (3) at (2,1)  {}; 
\node[vertex] (4) at (3,1)  {}; 
\node[vertex] (5) at (4,1)  {}; 
\node[vertex] (6) at (2,2)  {}; 

\draw (1) -- (2) -- (3) -- (4) -- (5);
\draw (2) -- (6) -- (4); \draw (3) -- (6);

\draw(2,.3) node {$F_6$};
\end{tikzpicture}
\qquad
\begin{tikzpicture}[scale=.55] %,thick]
%\tikzstyle{vertex}=[circle,inner sep=2pt,fill=black];
\tikzstyle{vertex}=[draw,circle,inner sep=1.5pt] 
\node[vertex] (1) at (0,1)  {};
\node[vertex] (2) at (1,1)  {};
\node[vertex] (3) at (2,1)  {}; 
\node[vertex] (4) at (3,1)  {}; 
\node[vertex] (5) at (4,1)  {}; 
\node[vertex] (6) at (5,1)  {};
\node[vertex] (7) at (2.5,2)  {};  

\draw (1) -- (2) -- (3) -- (4) -- (5) -- (6);
\draw (2) -- (7) -- (3); \draw (4) -- (7) -- (5);

\draw(2.5,.3) node {$F_7$};
\end{tikzpicture}
\qquad
\begin{tikzpicture}[scale=.55] %,thick]
%\tikzstyle{vertex}=[circle,inner sep=2pt,fill=black];
\tikzstyle{vertex}=[draw,circle,inner sep=1.5pt] 
\node[vertex] (1) at (0,1)  {};
\node[vertex] (2) at (1,1)  {};
\node[vertex] (3) at (2,1)  {}; 
\node[vertex] (4) at (3,1)  {}; 
\node[vertex] (5) at (0.5,2)  {}; 
\node[vertex] (6) at (1.5,2)  {};
\node[vertex] (7) at (2.5,2)  {};  

\draw (1) -- (2) -- (3) -- (4) -- (7) -- (6) -- (5) -- (1);
\draw (5) -- (2) -- (6) -- (3) -- (7);

\draw(1.5,.3) node {$F_8$};
\end{tikzpicture}
\qquad
\begin{tikzpicture}[scale=.55] %,thick]
%\tikzstyle{vertex}=[circle,inner sep=2pt,fill=black];
\tikzstyle{vertex}=[draw,circle,inner sep=1.5pt] 
\node[vertex] (1) at (0,1)  {};
\node[vertex] (2) at (2,1)  {};
\node[vertex] (3) at (4,1)  {}; 
\node[vertex] (7) at (2,2.3)  {}; 
\node[vertex] (4) at (1,3)  {}; 
\node[vertex] (5) at (3,3)  {};
\node[vertex] (6) at (5,3)  {};  

\draw (1) -- (2) -- (3) -- (6) -- (5) -- (4) -- (1);
\draw (4) -- (2) -- (5) -- (3); \draw (2) -- (7) -- (3); \draw (4) -- (7) -- (5);

\draw(2.5,.3) node {$F_9$};
\end{tikzpicture}
\end{center}
\caption{Forbidden induced subgraphs.}\label{fig:forbidden} 
\end{figure} 

Our main results are as follows.

\begin{thm}\label{thm:forest}
The Gallai graph $\Gamma(G)$ of a graph $G$ is a forest
if and only if $G$ is an $(F_1,\ldots, F_9)$-free chordal graph.
\end{thm} 
The {\it gem} is the graph that arises by removing the two end-vertices from $F_7$.
A set $U$ of vertices of a graph $G$ is {\it homogeneous}
if every vertex in $V(G)\setminus U$ is adjacent 
either to all vertices in $U$
or to no vertex in $U$.
A homogeneous set $U$ is {\it non-trivial} if $|U|\not\in \{ 0,1,|V(G)|\}$.

\begin{figure}[H]
\begin{center}
\begin{tikzpicture}[scale=.55] %,thick]
%\tikzstyle{vertex}=[circle,inner sep=2pt,fill=black];
\tikzstyle{vertex}=[draw,circle,inner sep=1.5pt] 
\node[vertex] (1) at (0,1)  {};
\draw(0,0.5) node {$a$};
\node[vertex] (2) at (1,1)  {};
\draw(1,0.5) node {$b$};
\node[vertex] (3) at (2,1)  {}; 
\draw(2,0.5) node {$e$};
\node[vertex] (4) at (2.5,2)  {}; 
\draw(2.5,2.5) node {$f$};
\node[vertex] (5) at (1.5,2)  {}; 
\draw(1.5,2.5) node {$c$};
\node[vertex] (6) at (0.5,2)  {};
\draw(0.5,2.5) node {$d$};

\draw (1) -- (2) -- (3) -- (4) -- (5) -- (6) -- (1);
\draw (6) -- (2) -- (5) -- (3);

\draw(1.25,-0.5) node {$F_8^-$};
\end{tikzpicture}
\end{center}
\caption{The graph $F_8^-$.}\label{fig:F8-} 
\end{figure}
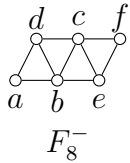 

\begin{thm}\label{thm:tree}
For a graph $G$ without isolated vertices, the following statements are equivalent:
\begin{itemize}
\item[\em (i)] The Gallai graph $\Gamma(G)$ of $G$ is a tree. 
\item[\em (ii)] Every non-trivial homogeneous set in $G$ is independent, 
and $G$ is an $(F_1,\ldots, F_9)$-free chordal graph.
\item[\em (iii)] $G$ is either the graph $F_8^-$ in Figure \ref{fig:F8-} 
or $G$ is connected and satisfies the following conditions:
\begin{itemize}
\item Every block of $G$ is isomorphic to $K_2$, $K_3$, or a gem.
\item Every cut-vertex of $G$ lies in at most two blocks and has degree at most $3$ in $G$.
\item Every block of $G$ that is isomorphic to $K_3$
has exactly two cut-vertices.
\item Every block of $G$ that is isomorphic to a gem
has exactly one cut-vertex.
\end{itemize}
\end{itemize}
\end{thm} 

The rest of the paper is devoted to the proofs of the above results.

\section{Proofs}

Before we proceed to the proofs of our results,
we collect some immediate observations.
\begin{itemize}
\item {\it Every graph is an induced subgraph of some Gallai graph.}
\end{itemize}
In fact, if $H$ is a graph 
and the graph $G$ has vertex set $V(H)\cup \{ x\}$ 
such that all vertices in $V(H)$ are neighbors of $x$ in $G$
and $G-x$ is the complement of $H$, 
then the subgraph of $\Gamma(G)$ induced by the edges of $G$ that are incident with $x$
is isomorphic to $H$.
This observation explains to some extend why the characterization of Gallai graphs is difficult.
\begin{itemize}
\item {\it If $G'$ is an induced subgraph of a graph $G$, 
then $\Gamma(G')$ is an induced subgraph of $\Gamma(G)$.}
\end{itemize}
This follows immediately from the definition.

For the convenience of the reader,
we include a proof of the following known result.

\begin{prop}[Le \cite{le1988}]\label{prop:fact2}
If $G$ is a graph without isolated vertices, then $\Gamma(G)$ is connected 
if and only if every non-trivial homogeneous set in $G$ is independent.
\end{prop}
{\it Proof:} Let $G$ be a graph without isolated vertices.

First we prove the necessity.
If $U$ is a non-trivial homogeneous set in $G$,
$uv$ is an edge of $G$ between two vertices in $U$,
and 
$xy$ is an edge of $G$ such that $x$ does not belong to $U$,
then $xy$ and $uv$ belong to distinct components of $\Gamma(G)$.
In fact, if $e_1\ldots e_\ell$ were a path in $\Gamma(G)$ with $e_1=xy$ and $e_\ell=uv$, 
then there is some index $i$  
such that $e_i$ joins a vertex $x_i$ in $V(G)\setminus U$ to a vertex $x_{i+1}$ in $U$
and $e_{i+1}$ joins $x_{i+1}$ to a vertex $x_{i+2}$ in $U$.
Since $U$ is homogeneous, $x_i$ is adjacent to $x_{i+2}$,
which implies the contradiction that $e_i$ and $e_{i+1}$ are not adjacent in $\Gamma(G)$.
This implies the necessity.

In order to prove the sufficiency, we assume that $\Gamma(G)$ is not connected.
Let $C$ be the vertex set of a component of $\Gamma(G)$,
that is, $C$ is a set of edges of $G$.
Let $V(C)$ denote the set of vertices of $G$ that are incident with an edge in $C$.
If $V(C)$ is a proper subset of $V(G)$,
then the definition of $V(C)$ implies that $V(C)$ is homogeneous,
that is, in this case $G$ has a non-trivial homogeneous set that is not independent.
Hence, we may assume that $V(C)=V(G)$ for all vertex sets $C$ of components of $\Gamma(G)$.
Now Lemma 4 in \cite{arju} implies a contradiction,
which completes the proof of the sufficiency. $\Box$

\medskip

\noindent We proceed to the proof of our first main result.

\medskip

\noindent {\it Proof of Theorem \ref{thm:forest}:}
Since the Gallai graph of a chordless cycle of length at least $4$
and of each of the graphs $F_1,\ldots,F_9$ contains a cycle,
the necessity follows.
In order to show the sufficiency,
let $G$ be an $(F_1,\ldots,F_9)$-free chordal graph.
We prove that $\Gamma(G)$ is a forest.
Clearly, we may assume that $G$ is connected.

\begin{claim}\label{claim:forest:1}
If $G$ contains an induced $F_8^-$,
then $G$ is isomorphic to $F_8^-$.
\end{claim}
{\it Proof of Claim \ref{claim:forest:1}:}
We denote the vertices of the induced $F_8^-$ as in Figure \ref{fig:F8-}.
For a contradiction, we assume that $G$ is not isomorphic to $F_8^-$.
Since $G$ is connected, some vertex $g$ in $V(G)\setminus V(F_8^-)$
is adjacent to some vertices in $V(F_8^-)$.
We consider different cases.

First we assume that $g$ is adjacent to $b$ but not to $a$.
Since $G[\{ a,b,c,g\}]$ and $G[\{ a,b,e,g\}]$ are no claws,
$g$ is adjacent to $c$ and $e$.
Since $G[\{ a,b,d,e,g\}]$ is not $F_2$,
$g$ is adjacent to $d$.
Since $G[V(F_8^-)\cup \{ g\}]$ is not $F_9$,
$g$ is adjacent to $f$.
Now $G-e$ is $F_3$,
which is a contradiction.

Next we assume that $g$ is adjacent to $a$ and $b$.
Since $G[\{ b,d,e,g\}]$ is not a claw,
$g$ is adjacent to $d$ or $e$.
Since $G[\{ a,c,d,g\}]$ and $G[\{ a,c,d,e,g\}]$ are no chordless cycles,
$g$ is adjacent to $d$.
Since $G[\{ a,b,c,e,g\}]$ is not $F_2$,
$g$ is adjacent to $c$ or $e$.
Since $G[\{ c,d,e,g\}]$ is not a chordless cycle,
$g$ is adjacent to $c$.
In view of the first case and the symmetry between $b$ and $c$,
we may assume that $g$ is adjacent to $f$.
Since $G[\{ b,e,f,g\}]$ is not a chordless cycle,
$g$ is adjacent to $e$.
Now $G[\{ a,d,e,f,g\}]$ is $F_2$,
which is a contradiction.

In view of the first two cases, 
we may assume that $g$ is not adjacent to $b$ and $c$.
If $g$ is adjacent to $d$ but not to $a$,
then $G[\{ a,d,c,g\}]$ is a claw.
Since $g$ has a neighbor in $V(F_8^-)$
and in view of the symmetry between $a$ and $f$,
we may assume that $g$ is adjacent to $a$.
Since $G[\{ a,b,e,g\}]$ is not a chordless cycle,
$g$ is not adjacent to $e$.
Since $G[\{ a,b,c,f,g\}]$ is not a chordless cycle,
$g$ is not adjacent to $f$.
Since $G[\{ a,b,c,d,f,g\}]$ is not $F_6$,
$g$ is adjacent to $d$.
Now $G[V(F_8^-)\cup \{ g\}]$ is $F_8$,
which is a contradiction.

This completes the proof of the claim.
$\Box$

\medskip

\noindent In the following 
we may assume that $G$ is $F_8^-$-free. 
We proceed by induction on the order of $G$.
Since the result holds for graphs of order at most $3$,
we assume that $G$ has order at least $4$.
We consider different cases.

\medskip

\noindent {\bf Case 1} {\it $G$ has an induced gem.}

\medskip

\noindent Let $a_1, \ldots, a_5$ be the vertices of an induced gem in $G$
such that $a_1a_2a_3a_4$ is the induced path of order $4$ in that gem. 
Let $A$ be the set of all vertices that are adjacent to $a_1$, $a_2$, $a_3$, and $a_4$. 
Since $a_5\in A$ and $G$ is chordal,
$A$ is a non-empty clique.

If a vertex $b$ not in $A\cup\{a_1,a_2,a_3,a_4\}$
has a neighbor $a$ in $A$,
then, since $G[\{ a,b,a_1,a_4\}]$ is not a claw,
we may assume, by symmetry, 
that $b$ is adjacent to $a_1$.
Since $G$ is chordal and $b$ is not adjacent to $a_1$, $a_2$, $a_3$, or $a_4$,
the vertex $b$ is not adjacent to $a_4$.
Since $G[\{ a,b,a_2,a_4\}]$ is not a claw,
$b$ is adjacent to $a_2$.
Since $G[\{ a,b,a_1,a_3,a_4\}]$ is not $F_2$,
$b$ is adjacent to $a_3$.
By symmetry, 
it follows that every vertex not in $A\cup\{a_1,a_2,a_3,a_4\}$ 
that has a neighbor in $A$ is
\begin{itemize}
\item either adjacent to $a_1,a_2,a_3$ and not adjacent to $a_4$ ({\it type 1})
\item or adjacent to $a_2,a_3,a_4$ and not adjacent to $a_1$ ({\it type 2}).
\end{itemize}
If $b$ is of type 1 and $b'$ is of type 2,
then $G[\{ a_1,a_2,a_3,a_4,b,b'\}]$ is either $F_3$ or $F_8^-$.
Hence, we assume that there is no vertex of type 2.
Let $B$ denote the set of vertices of type 1.
Since $G$ is chordal,
$A\cup B$ is a clique.  
Let $C=V(G)\setminus (A\cup B\cup\{a_1,a_2,a_3,a_4\})$.

If a vertex $v$ in $C$ has a neighbor $b$ in $B$,
then, since $G[\{ v,b,a_1,a_3\}]$ is not a claw,
$v$ is adjacent to $a_1$ or $a_3$.
Since $v$ is not adjacent to $a_5$ and $G[\{ a_5,v,a_1,a_3\}]$ is not a chordless cycle,
$v$ is not adjacent to $a_1$ or $a_3$.
Since $G[\{ a_4,a_5,b,v\}]$ is not a chordless cycle, 
$v$ is not adjacent to $a_4$.
If $v$ is adjacent to $a_1$ and not to $a_3$,
then $G[\{ v,b,a_1,a_3,a_4,a_5\}]$ is $F_8^-$.
If $v$ is adjacent to $a_3$ and not to $a_1$,
then $G[\{ v,b,a_1,a_3,a_4,a_5\}]$ is $F_5$.
This implies that no vertex in $C$
that has a neighbor in $B$.

If a vertex $v$ in $C$
that is adjacent to $a_2$ and $a_3$, 
then, since $G[\{ v,a_1,a_2,a_3,a_4,a_5\}]$ is not $F_5$,
$v$ is adjacent to $a_1$ or $a_4$,
which implies that either
$G[\{ a_1,v,a_3,a_5\}]$
or
$G[\{ a_2,v,a_4,a_5\}]$
is a chordless cycle.
This implies that no vertex in $C$
is adjacent to $a_2$ and $a_3$.

If a vertex $v$ in $C$
that is adjacent to $a_2$ and not adjacent to $a_3$, 
then, since $G[\{ v,a_1,a_2,a_3\}]$ is not a claw,
$v$ is adjacent to $a_1$.
Since $G$ is chordal,
$v$ is not adjacent to $a_4$.
Now $G[\{ v,a_1,a_2,a_3,a_4,a_5\}]$ is $F_8^-$.
This implies that no vertex in $C$
is adjacent to $a_2$ and not adjacent to $a_3$.
Similarly, it follows that
no vertex in $C$
is adjacent to $a_3$ and not adjacent to $a_2$.
Altogether, since $G$ is chordal,
the neighborhood of every vertex in $C$
in $A\cup B\cup\{a_1,a_2,a_3,a_4\}$
is either empty or $\{ a_1\}$ or $\{ a_4\}$.

If a vertex in $C$ is adjacent to $a_1$ and 
another vertex in $C$ is adjacent to $a_4$,
then $G$ contains $F_7$.
If $B$ is not empty and a vertex in $C$ is adjacent to $a_1$,
then $G$ contains $F_6$.
Hence, by symmetry, 
we may assume 
that no vertex in $C$ is adjacent to $a_1$.

By induction, $\Gamma(G-a_1)$ is a forest.
In view of the above observations, 
the edges of $G$ incident with $a_1$ form an independent set $X$ of $\Gamma(G)$,
and for every vertex $a$ in $A$, 
the edge $aa_3$ is an isolated vertex of $\Gamma(G-a_1)$.
Since $\Gamma(G)$ arises from the disjoint union of $\Gamma(G-a_1)$ and $X$ by 
\begin{itemize}
\item adding the two edges $(a_1a)(aa_3)$ and $(a_1a)(aa_4)$ for every $a\in A$,
\item adding the edge $(a_1b)(ba_3)$ for every $b\in B$, and
\item adding the edge $(a_1a_2)(a_2a_3)$,
\end{itemize}
$\Gamma(G)$ is a forest,
which completes the proof in Case 1.

\medskip

\noindent In view of Case 1 
we may now assume that $G$ has no induced gem.
Two distinct vertices $x$ and $y$ of $G$ with $N_G(x)\setminus \{ y\}=N_G(y)\setminus \{ x\}$
are called {\it twins}.

\medskip

\noindent {\bf Case 2} {\it $G$ contains two distinct vertices $x$ and $y$ that are twins.}

\medskip

\noindent Let $C=N_G(x)\setminus \{ y\}$.
If $x$ and $y$ are not adjacent, 
then, since $G$ is chordal, $C$ is a clique.
Since $G$ is claw-free, we obtain $V(G)=C\cup \{ x,y\}$, 
and $\Gamma(G)$ is a forest
that consists of ${|C|\choose 2}$ isolated vertices
and $|C|$ components of order $2$.
Hence, we may assume that $x$ and $y$ are adjacent.

If $c_1c_2\ldots c_{\ell}c_1$ is a chordless cycle in $\overline{G}[C]$,
then, since $G$ is claw-free, $\ell\geq 4$.
If $\ell=4$, then $G[ \{ x\}\cup C]$ is $F_2$,
if $\ell=5$, then $G[C]$ is $C_5$, and
if $\ell\geq 6$, then $G[ \{ c_1,c_5,c_2,c_4\}]$ is a chordless cycle.
Altogether, it follows that $\overline{G}[C]$ is a forest.

Let $A\subseteq C$ be the set of vertices in $C$ that have a neighbor not in $\{ x,y\}\cup C$.
Let $B=C\setminus A$.
If $A$ is empty, then $\Gamma(G)$ is the disjoint union
of an isolated vertex $xy$,
two disjoint copies of the forest $\overline{G}[C]$ induced by the edges of $G$ joining $\{ x,y\}$ to $C$,
and the graph $\Gamma(G-\{ x,y\})$, which is a forest by induction,
that is, $\Gamma(G)$ is a forest.
Hence, we may assume that $A$ is not empty.

Let $D$ be the set of vertices not in $\{ x,y\}\cup C$ that have a neighbor in $A$.
By definition, every vertex in $A$ has a neighbor in $D$.

If a vertex $a$ in $A$ has two neighbors $d_1$ and $d_2$ in $D$, 
then 
either $d_1$ and $d_2$ are not adjacent and $G[\{ a,d_1,d_2,x\}]$ is a claw,
or $d_1$ and $d_2$ are adjacent and $G[\{ a,d_1,d_2,x,y\}]$ is $F_2$.
Hence every vertex in $A$ has exactly one neighbor in $D$.

If two distinct vertices $a_1$ and $a_2$ in $A$ are not adjacent, 
then let $a_1'$ and $a_2'$ denote their neighbors in $D$, respectively.
If $a_1'$ and $a_2'$ are equal or adjacent, then $G[\{ x,a_1,a_2,a_1',a_2'\}]$ is a chordless cycle.
If $a_1'$ and $a_2'$ are distinct and not adjacent, then $G[\{ x,y,a_1,a_2,a_1',a_2'\}]$ is $F_6$.
This implies that $A$ is a clique.

If two vertices $a_1'$ and $a_2'$ in $D$ are adjacent, 
then let $a_1$ and $a_2$ denote their neighbors in $A$, respectively.
The graph $G[\{ a_1,a_2,a_1',a_2'\}]$ is a chordless cycle.
This implies that $D$ is an independent set.

First we assume that $A$ has only one element $a$.
Let $a'$ be the unique element of $D$.
By induction, the graph $\Gamma(G-x)$ is a forest.
Since $\Gamma(G)$ arises from $\Gamma(G-x)$
by adding 
\begin{itemize}
\item the edge $xy$ as an isolated vertex,
\item a disjoint copy of the subforest of $\Gamma(G-x)$ 
induced by the edges of $G$ incident with $y$, and 
\item an edge between $xa$ and $aa'$,
\end{itemize}
$\Gamma(G)$ is a forest.
Hence, we may assume that $A$ has at least two elements.

If some vertex $b$ in $B$ is not adjacent to some vertex $a_1$ in $A$,
then let $a_2$ be a vertex in $A$ distinct from $a_1$.
Let $a_1'$ and $a_2'$ denote the neighbors of $a_1$ and $a_2$ in $D$, respectively.
If $b$ is not adjacent to $a_2$, 
then 
either $a_1'=a_2'$ and $G[\{ x,y,a_1,a_2,a_1',b\}]$ is $F_3$
or 
$a_1'\not=a_2'$ and $G[\{ x,a_1,a_2,a_1',a_2',b\}]$ is $F_4$.
Hence $b$ is adjacent to $a_2$.
Now
either $a_1'=a_2'$ and $G[\{ x,a_1',a_1,a_2,b\}]$ is a gem 
or $a_1'\not=a_2'$ and $G[\{ a_1,a_2,a_2',b\}]$ is a claw.
This implies that every vertex in $B$ is adjacent to every vertex in $A$.

If two distinct vertices $b$ and $b'$ in $B$ are not adjacent,
then $b$, $b'$, a vertex in $A$ and its neighbor in $D$ induce a claw in $G$.
Hence $C$ is a clique.
By induction, the graph $\Gamma(G-x)$ is a forest.
Since $\Gamma(G)$ arises from $\Gamma(G-x)$ by adding
\begin{itemize}
\item the edge $xy$ as an isolated vertex,
\item for each $b\in B$, the edge $xb$ as an isolated vertex, and
\item for each $a\in A$ whose neighbor in $D$ is $a'$,
the edge $xa$ as an end-vertex that is adjacent only to $aa'$, 
\end{itemize}
$\Gamma(G)$ is a forest.
This completes the proof in Case 2.

\medskip

\noindent In view of Cases 1 and 2, 
we may assume that 
$G$ is a gem-free twin-free chordal graph.
By a result of Howorka \cite{howo}, 
$G$ is distance-hereditary
and, by a result of Bandelt and Mulder \cite{bamu},
$G$ has a vertex of degree $1$,
which leads us to our final case.

\medskip

\noindent {\bf Case 3} {\it $G$ contains a vertex $v$ of degree $1$.}

\medskip

\noindent Let $w$ denote the neighbor of $v$.
Let $Q=N_G(w)\setminus \{ v\}$.
If $Q$ has just one element $q$, 
then $\Gamma(G)$ arises from $\Gamma(G-v)$
by adding the vertex $vw$ and an edge between $vw$ and $wq$.
Since $\Gamma(G-v)$ is a forest by induction, also $\Gamma(G)$ is a forest.
Hence, we may assume that $Q$ has at least two elements.
Since $G$ is claw-free, $Q$ is a clique.

If a vertex $z$ not in $N_G[w]$ has two neighbors $q_1$ and $q_2$ in $Q$,
then, since $G$ is twin-free, 
we may assume that $q_1$ has a neighbor $q_1'$ that is not adjacent to $q_2$.
If $q_1'$ and $z$ are not adjacent, then $G[\{ q_1,w,z,q_1'\}]$ is a claw,
and 
if $q_1'$ and $z$ are adjacent, then $G[\{ q_1,q_2,w,z,q_1'\}]$ is a gem.
This implies that every vertex not in $N_G[w]$
has at most one neighbor in $Q$.

If no vertex in $Q$ has a neighbor not in $N_G[w]$,
then $\Gamma(G)$ is the disjoint union of ${|Q|\choose 2}$ isolated vertices
and a star of order $|Q|+1$.
Hence, we may assume that some vertex in $Q$ has a neighbor not in $N_G[w]$.

If two vertices $q_1$ and $q_2$ in $Q$ have neighbors, say $q_1'$ and $q_2'$, respectively, not in $N_G[w]$,
then either $q_1'$ and $q_2'$ are adjacent and $G[\{ q_1,q_2,q_1',q_2'\}]$ is a chordless cycle
or $q_1'$ and $q_2'$ are not adjacent and $G[\{ v,w,q_1,q_2,q_1',q_2'\}]$ is $F_4$.
Hence exactly one vertex in $Q$, say $q_1$, has a neighbor not in $N_G[w]$.
Since $G$ is twin-free, $Q$ contains exactly one further element $q_2$.

Since $\Gamma(G-v)$ is a forest by induction, 
$wq_2$ is an isolated vertex in $\Gamma(G-v)$, and
$\Gamma(G)$ arises from $\Gamma(G-v)$
by adding the vertex $vw$ and the two edges $(vw)(wq_1)$ and $(vw)(wq_2)$,
we obtain that $\Gamma(G)$ is a forest.
This completes the proof in Case 3,
which completes the entire proof. $\Box$

\medskip

\noindent We proceed to the proof of our second main result.

Clearly, 
Theorem \ref{thm:forest} and Proposition \ref{prop:fact2} imply the equivalence of 
(i) and (ii) in Theorem \ref{thm:tree}.
In order to complete the proof of Theorem \ref{thm:tree} 
it would suffice 
to prove the equivalence of (ii) and (iii) in Theorem \ref{thm:tree}.
Since we want to emphasize the interplay between potential cycles in $\Gamma(G)$ 
and the structural features expressed in (iii), 
we complete the proof of Theorem \ref{thm:tree} 
by showing the equivalence of (i) and (iii) in Theorem \ref{thm:tree} directly,
which might be slightly longer yet more instructive.

\medskip

\noindent {\it Proof of Theorem \ref{thm:tree}:}
Theorem \ref{thm:forest} and Proposition \ref{prop:fact2} imply the equivalence of (i) and (ii). 
We proceed to the proof of the equivalence of (i) and (iii).

We first prove the sufficiency, that is, that (iii) implies (i). 
Since the Gallai graph of the graph in Figure \ref{fig:F8-} is a tree, 
we may assume that $G$ is connected satisfies the four conditions stated in (iii).
Removing all vertices of $G$ that are no cut-vertices results in a path $P$
whose Gallai graph $\Gamma(P)$ is again a path.
If $B$ is a block of $G$ that is isomorphic to $K_3$,
then let $u$ and $v$ denote the two cut-vertices of $G$ in $B$
and let $w$ denote the third vertex of $B$.
The conditions imply that $u$ has exactly one neighbor $u'$ that is not in $B$
and that $v$ has exactly one neighbor $v'$ that is not in $B$.
In $\Gamma(G)$ the two edges $uw$ and $vw$ are end-vertices 
adjacent to $uu'$ and $vv'$, respectively.
If $B$ is a block of $G$ that is isomorphic to a gem,
then let $u$ denote the unique cut-vertex of $G$ in $B$.
The conditions imply that $u$ has degree $2$ in $B$
and has exactly one neighbor $u'$ that is not in $B$.
The two vertices $u$ and $u'$ form a block of $G$
and the edge $uu'$ is an end-vertex of $\Gamma(P)$.
In $\Gamma(G)$ the edges of $B$ form two small subtrees 
each attached by one edge to the end-vertex $uu'$ of $\Gamma(P)$;
one is isomorphic to $P_3$ and the other one is isomorphic to $P_4$.
Altogether, $\Gamma(G)$ is a tree.

Now we prove the necessity.
Therefore, let $G$ be a graph without an isolated vertex
such that $\Gamma(G)$ is a tree.
Clearly, we may assume that $G$ is not the graph in Figure \ref{fig:F8-}.
Since the Gallai graph of a chordless cycle of length at least $4$ is a cycle,
the graph $G$ is chordal.

We begin with a useful observation.

\setcounter{claim}{0}

\begin{claim}\label{claim1}
If $H$ is a proper induced subgraph of $G$ without an isolated vertex
such that $\Gamma(H)$ is connected, 
then there is a vertex in $V(G)\setminus V(H)$ 
that is adjacent to some vertices of $H$ but not to all.
Furthermore, 
if $u$ in $V(G)\setminus V(H)$ is adjacent to $v$ in $V(H)$,
then $u$ is adjacent to all but at most one neighbor of $v$ in $H$.
\end{claim}
{\it Proof of Claim \ref{claim1}:}
The Gallai graph $\Gamma(H)$ of $H$ is a subtree of $\Gamma(G)$.
Since $H$ does not contain all edges of $G$ and $\Gamma(G)$ is connected,
some edge of $H$ spans a $P_3$ with an edge joining 
$V(G)\setminus V(H)$ to $V(H)$,
which implies the first part of the claim.
If $u$ in $V(G)\setminus V(H)$ is adjacent to $v$ in $V(H)$
but is not adjacent to two neighbors, say $x$ and $y$, of $v$ in $H$,
then $uv$ is adjacent in $\Gamma(G)$ to $vx$ and $vy$.
Since $\Gamma(H)$ contains a path between $vx$ and $vy$,
$\Gamma(G)$ contains a chordless cycle,
which is a contadiction and completes the proof of Claim \ref{claim1}.
$\Box$

\medskip

\noindent It is a trivial consequence of Claim \ref{claim1}
that $G$ is claw-free.

Our next claim concerns induced diamonds in $G$.

\begin{claim}\label{claim1.5}
If $G$ contains an induced diamond $D$,
then there is a way of naming the vertices of $D$ as $a$, $b$, $c$, and $d$
such that $E(D)=\{ ab,bc,cd,da,bd\}$
and there are two vertices $e$ and $f$ in $G$ with 
$N_G(e)\cap \{ a,b,c,d\}=\{ b,c\}$
and 
$N_G(f)\cap \{ a,b,c,d,e\}=\{ e\}$ such that 
every vertex $g$ in $V(G)\setminus \{ a,b,c,d,e,f\}$
that is adjacent to a vertex in $\{ a,b,c,d,e,f\}$ is adjacent exactly to $f$.
\end{claim}
{\it Proof of Claim \ref{claim1.5}:}
Let $G$ contain an induced diamond
with vertices $a$, $b$, $c$ and $d$,
and edges $ab$, $bc$, $cd$, $da$, and $bd$.
Claim \ref{claim1} applied to $G[\{ b,d\}]$ 
implies the existence of a vertex $e$ 
that is adjacent to exactly one of the two vertices $b$ and $d$.
By symmetry, we may assume that $e$ is adjacent to $b$ but not to $d$.
Since $G[\{ a,b,c,e\}]$ is not a claw, $e$ is adjacent to $a$ or $c$.
Since $G[\{ a,c,d,e\}]$ is not a chordless cycle, $e$ is not adjacent to both $a$ and $c$.
By symmetry, we may assume that $e$ is adjacent to $c$ but not to $a$,
that is, $G[\{ a,b,c,d,e\}]$ is a gem.

To complete the proof of Claim \ref{claim1.5},
we establish a further claim.

\begin{claim}\label{claim2}
There is a vertex that is adjacent to some but not all vertices in $\{ a,c,d,e\}$
and is not adjacent to $b$.
\end{claim}
{\it Proof of Claim \ref{claim2}:}
We call a vertex that is adjacent to some but not all vertices in $\{ a,c,d,e\}$ {\it interesting}.
By Claim \ref{claim1} applied to $G[\{ a,c,d,e\}]$, there is at least one interesting vertex.
For a contradiction, we assume that every interesting vertex is adjacent to $b$.
We first show that every interesting vertex is adjacent to $c$ and $d$.
For a contradiction, we may assume, by symmetry, 
that the interesting vertex $f$ is not adjacent to $c$.
Since $G[\{ a,b,c,f\}]$ is not a claw, $f$ is adjacent to $a$.
By Claim \ref{claim1} applied to $G[\{ a,b,c,e\}]$, 
$f$ is adjacent to $e$.
Since $adcefa$ is not a chordless cycle in $G$,
$f$ is adjacent to $d$.
Now $dcefd$ is a chordless cycle in $G$,
which is a contradiction,
that is, every interesting vertex is adjacent to $c$ and to $d$.
Since for an interesting vertex $f$,
$G[\{ a,b,e,f\}]$ is not a claw,
every interesting vertex is adjacent to $a$ or $e$,
that is, for every interesting vertex $f$,
the set $N_G(f)\cap \{ a,c,d,e\}$ is either $\{ c,d,e\}$ or $\{ a,c,d\}$.

Let $f$ be an interesting vertex.
By symmetry, we may assume that $N_G(f)\cap \{ a,c,d,e\}=\{ c,d,e\}$.

By Claim \ref{claim1} applied to $G[\{ c,f\}]$, 
there is a vertex $g$ distinct from $c$ and $f$ 
that is adjacent to exactly one of $c$ and $f$.
First, we assume that $g$ is adjacent to $f$ but not to $c$.
Since $G[\{ d,e,f,g\}]$ is not a claw, $g$ is adjacent to $d$ or $e$.
Hence $g$ is interesting, which implies a contradiction 
to the above observation that every interesting vertex is adjacent to $c$.
Hence $g$ is adjacent to $c$ but not to $f$.
If $g$ is not interesting, that is,
$g$ is adjacent to all vertices in $\{a,c,d,e\}$,
then, since $abega$ is not a chordless cycle in $G$, 
$g$ is adjacent to $b$.
Now $(eg)(ef)(df)(dg)(eg)$ is a cycle in $\Gamma(G)$,
which is a contradiction.
Hence $g$ is interesting and thus adjacent to $b$.
Since $G[\{ a,b,f,g\}]$ is not a claw, 
$g$ is adjacent to $a$.
Since $g$ is interesting, 
this implies that $N_G(g)\cap \{ a,c,d,e\}=\{ a,c,d\}$.
Now $(bg)(be)(ab)(bf)(bg)$ is a cycle in $\Gamma(G)$,
which is a contradiction
and completes the proof of Claim \ref{claim2}. $\Box$
 
\medskip

\noindent Let $f$ the vertex whose existence is guaranteed by Claim \ref{claim2}.
Let $H=G[\{ a,b,c,d,e,f\}]$.

First we assume that $f$ is adjacent to $c$.
Since $abcfa$ is not a chordless cycle in $G$, $f$ is not adjacent to $a$.
If $f$ is adjacent to $d$, then, 
since $dbefd$ is not a chordless cycle in $G$, $f$ is not adjacent to $e$
and $(ad)(df)(bd)(be)(ab)(bc)(cf)(ce)(cd)(ad)$ is a cycle in $\Gamma(G)$,
which is a contradiction.
Hence $f$ is not adjacent to $d$.
Since $G[\{ c,d,e,f\}]$ is not a claw, $f$ is adjacent to $e$.
Since $H$ is the graph in Figure \ref{fig:F8-},
$H$ is a proper induced subgraph of $G$.
Let $g$ be the vertex whose existence is guaranteed 
by Claim \ref{claim1} applied to $H$.
By the second half of Claim \ref{claim1},
that $g$ is adjacent to $b$ or $c$.
By symmetry, we assume that $g$ is adjacent to $c$.
Iteratively applying the second half of Claim \ref{claim1},
it follows in turn that $g$ is adjacent to $b$, $d$, and $e$.
Hence 
$N_G(g)\cap V(H)\in \{ V(H)\setminus \{ a\},V(H)\setminus \{ f\},V(H)\setminus \{ a,f\}\}$.
If $N_G(g)\cap V(H)=V(H)\setminus \{ f\}$, 
then $(eg)(ef)(be)(ab)(bc)(cf)(cg)(ag)(eg)$ is a cycle in $\Gamma(H)$,
which is a contradiction.
If $N_G(g)\cap V(H)=V(H)\setminus \{ a\}$, 
a symmetric argument leads to a contradiction.
If $N_G(g)\cap V(H)=V(H)\setminus \{ a,f\}$, 
then $(eg)(ef)(be)(ab)(bc)(cf)(cd)(ad)(dg)(eg)$ is a cycle in $\Gamma(H)$,
which is a contradiction.
Hence, by symmetry, we may assume that $f$ is not adjacent to $c$ and $d$.
By symmetry, we may further assume that $f$ is adjacent to $e$.
Since $abefa$ is not a chordless cycle in $G$, $f$ is not adjacent to $a$,
that is, $N_G(f)\cap \{ a,b,c,d,e\}=\{ e\}$.
Note that the vertices $e$ and $f$ are as required in the statement of Claim \ref{claim1.5}
and that $\Gamma(H)$ is a tree.

Let $g$ be a vertex $V(G)\setminus V(H)$ that is adjacent to a vertex in $V(H)$.
For a contradiction, we assume that $N_G(g)\cap V(H)\not=\{ f\}$,
that is, $g$ has a neighbor in $\{ a,b,c,d,e\}$.
Since $G[\{ a,c,f,g\}]$ is not a claw, $g$ is not adjacent to at least one vertex in $\{ a,c,f\}$.
If $g$ is not adjacent to $c$,
then iteratively applying the second half of Claim \ref{claim1} to $H$
implies that $N_G(g)\cap V(H)=V(H)\setminus \{ c\}$ and
the path $(bc)(bg)(fg)(dg)(cd)$ in $\Gamma(G)$ 
together with a path in $\Gamma(H)$ between $bc$ and $cd$
forms a cycle in $\Gamma(G)$, which is a contradiction.
Hence $g$ is adjacent to $c$.
If $g$ is not adjacent to $a$,
then iteratively applying the second half of Claim \ref{claim1} to $H$
implies that $N_G(g)\cap V(H)\in \{ V(H)\setminus \{ a\},V(H)\setminus \{ a,f\}\}$.
If $N_G(g)\cap V(H)=V(H)\setminus \{ a\}$,
the path $(ab)(bg)(fg)(dg)(ad)$ in $\Gamma(G)$ 
together with a path in $\Gamma(H)$ between $ab$ and $ad$
forms a cycle in $\Gamma(G)$, which is a contradiction.
If $N_G(g)\cap V(H)=V(H)\setminus \{ a,f\}$,
the path $(ad)(dg)(eg)(ef)$ in $\Gamma(G)$ 
together with a path in $\Gamma(H)$ between $ad$ and $ef$
forms a cycle in $\Gamma(G)$, which is a contradiction.
Hence $g$ is adjacent to $a$ and not adjacent to $f$. 
If $g$ is adjacent to $e$ but not adjacent to $d$, 
then the path $(ef)(eg)(ag)(ad)$ in $\Gamma(G)$ 
together with a path in $\Gamma(H)$ between $ef$ and $ad$
forms a cycle in $\Gamma(G)$, which is a contradiction.
This together with Claim \ref{claim1} applied to $H$ implies that 
$N_G(g)\cap V(H)\in \{ V(H)\setminus \{ f\},V(H)\setminus \{ e,f\}\}$.
If $N_G(g)\cap V(H)=V(H)\setminus \{ e,f\}$, 
then $\Gamma(G[\{ a,b,c,d,e,f,g\}])$ is a forest with two components 
one of which is the isolated vertex $dg$.
By Claim \ref{claim1} applied to $G[\{ d,g\}]$,
there is a vertex $h$ that is adjacent to exactly one of $d$ and $g$.
If $h$ is adjacent to $g$ but not to $d$, 
then $h$ is adjacent to one of $b$ and $c$,
that is, $h$ is a vertex that is adjacent to a vertex in $V(H)$
but satisfies $N_G(h)\cap V(H) \not=\{ f\}$.
Applying the same arguments to $h$ 
as we applied above to $g$ yields a contradiction.
Hence $h$ is adjacent to $d$ but not to $g$.
By symmetry with $g$, we obtain
$N_G(h)\cap V(H)\in \{ V(H)\setminus \{ f\},V(H)\setminus \{ e,f\}\}$.
Since $G[\{ c,e,g,h\}]$ is not a claw, 
we have 
$N_G(h)\cap V(H)=V(H)\setminus \{ f\}$.
Now the path $(ef)(eh)(ah)(ag)(cg)(ce)$ in $\Gamma(G)$ 
together with a path in $\Gamma(H)$ between $ef$ and $ce$
forms a cycle in $\Gamma(G)$, which is a contradiction.
Altogether we obtain that $N_G(g)\cap V(H)=V(H)\setminus \{ f\}$.
By symmetry, this easily implies that $g$ is a true twin of $b$,
which implies the contradiction, 
that $bg$ is isolated in $\Gamma(G)$.
This completes the proof of Claim \ref{claim1.5}. $\Box$

\medskip

\noindent Let $G'$ arise from $G$ by deleting all vertices of $G$
that belong to an induced diamond of $G$ 
that contains a vertex of degree $2$.
By Claim \ref{claim1.5}, 
$G'$ is a diamond-free graph 
such that $\Gamma(G')$ is a subtree of $\Gamma(G)$
and in order to complete the proof, 
it suffices to show that the conditions stated in the theorem hold for $G'$.
We establish some properties of the blocks of $G'$.

Since every $2$-connected diamond-free chordal graph is complete,
all blocks of $G'$ are complete.
If some block $B$ of $G'$ contains $3$ cut-vertices $a$, $b$, and $c$
and $a'$, $b'$, and $c'$ are neighbors of $a$, $b$, and $c$ outside of $B$,
respectively, then $(aa')(ab)(bb')(bc)(cc')(ac)(aa')$ is a cycle in $\Gamma(G')$,
which is a contradiction.
Hence every block of $G$ has at most $2$ cut-vertices.
Since $G'$ does not contain true twins, 
this implies that all blocks of $G'$ are isomorphic to $K_2$ or $K_3$.
Since $G'$ is claw-free, 
every cut-vertex of $G'$ lies in at most two blocks.
If some cut-vertex of $G'$ has degree at least $4$,
then two blocks of $G'$ of order $3$ share a cut-vertex 
and $\Gamma(G')$ contains a cycle of length $4$,
which is a contradiction.
Hence every cut-vertex of $G'$ has degree at most $3$,
which completes the proof of the theorem. $\Box$

\medskip

\noindent The main open problem related to Gallai graphs 
is the complexity of their recognition
and/or their efficient characterization.

\end{document}